\documentclass[12pt]{article}
\usepackage{epsfig,latexsym,amsfonts,amsmath,geometry}
\usepackage{colortbl}%

\def\R{\mathbb{R}}

\def\f{\varphi}

\def\irn{\int\limits_{\R^n}}

\def\Dtau{\left(-\Delta_{\R^n}\right)^{\!\tau}\!}

\def\Ds{\left(-\Delta_{\R^n}\!\right)^{\!s}\!}

\def\Dshalf{\left(-\Delta_{\R^n}\!\right)^{\!\frac{s}{2}}\!}

\def\eps{\varepsilon}

\def\iirn{\iint\limits_{\R^n\!\times\R^n}}

\def\irn{\int\limits_{\R^n}}

\def\proof{\noindent{\textbf{Proof. }}}
\def\QED{\hfill {$\square$}\goodbreak \medskip}

\newtheorem{Theorem}{Theorem}
\newtheorem{Lemma}
{Lemma}
\newtheorem{Proposition}
{Proposition}
\newtheorem{Remark}
{Remark}
\newtheorem{Example}
{Example}

\begin{document}

\title 
{A note on truncations in fractional Sobolev spaces}

\author{Roberta Musina\footnote{Dipartimento di Scienze Matematiche, Informatiche e Fisiche, Universit\`a di Udine,
via delle Scienze, 206 -- 33100 Udine, Italy. Email: {roberta.musina@uniud.it}. 
{Partially supported by Miur-PRIN 2015233N54.}}~ and
Alexander I. Nazarov\footnote{St.Petersburg Department of Steklov Institute, Fontanka, 27, St.Petersburg, 191023, Russia
and St.Petersburg State University, 
Universitetskii pr. 28, St.Petersburg, 198504, Russia. E-mail: {al.il.nazarov@gmail.com}. }
}
\date{}

\maketitle

\begin{abstract}
\footnotesize
\noindent
We study the Nemytskii operators $u\mapsto |u|$ and $u\mapsto u^{\pm}$ 
in  fractional Sobolev spaces $H^s(\R^n)$, $s>1$.

\bigskip

\noindent
\textbf{Keywords:} {Fractional Laplacian - Sobolev spaces - Truncation operators}

\bigskip

\noindent
\textit{2010 Mathematics Subject Classification:} {46E35, 47H30.}
\end{abstract}

\section{Introduction. Main result}

In this paper we discuss the relation between the map $u\mapsto |u|$ and the {\em Dirichlet Laplacian}. Recall that
the Dirichlet Laplacian $\Ds u$ of order $s>0$ of a function $u\in L^2(\R^n)$, $n\ge 1$,  is the distribution 
$$
\langle \Ds u,\f\rangle\equiv\irn u~\!\Ds\f~\!dx:=\irn|\xi|^{2s}{\mathcal F}[\f]~\!\overline{{\mathcal F}[u]}~\!d\xi~,\quad 
\f\in {\cal C}^\infty_0(\R^n),
$$
where $$
{\mathcal F}[u](\xi)= (2\pi)^{-\frac{n}{2}}
\irn e^{-i~\!\!\xi\cdot x}u(x)~\!dx
$$ 
is the Fourier transform in $\R^n$. The Sobolev--Slobodetskii space
$$
 H^s(\R^n) =\{u\in L^2(\R^n)~|~\Dshalf u\in L^2(\R^n)~\}
 $$
naturally inherits an Hilbertian structure from the scalar product
$$
 (u,v)=\langle\Ds u,v\rangle+\irn uv~\!dx~\!.
$$
The standard reference for the operator $\Ds$ and  functions in $H^s(\R^n)$ is the monograph
\cite{Tr} by Triebel. 

For any positive order $s\notin \mathbb N$ we introduce the constant 
\begin{equation}
\label{eq:Cns}
C_{n,s}=\frac{2^{2s}s}{\pi^{\frac{n}{2}}} ~\!\frac{\Gamma\big(\frac{n}{2}+s\big)}{\Gamma\big(1-s\big)}~\!.
\end{equation}
Notice that 
\begin{equation}
\label{eq:sign}
C_{n,s}>0\quad \text{if  ~$\lfloor s\rfloor$~ is even}; \qquad C_{n,s}<0\quad \text{if  ~$\lfloor s\rfloor$~ is odd,}
\end{equation}
where $ \lfloor s\rfloor$ stands for the integer part of $s$. It is well known that for $s\in(0,1)$
and $u,v\in H^s(\R^n)$ one has
\begin{equation}
\label{eq:hypersing}
\langle \Ds  u,v\rangle= \frac{C_{n,s}}{2}\irn\frac{(u(x)-u(y))(v(x)-v(y))}{|x-y|^{n+2s}}~\!dx~\!.
\end{equation}

Let us recall some known facts about the Nemytskii operator $|\cdot|: u\mapsto |u|$. 

\medskip
\noindent
1. 
 $|\cdot |$ is a Lipschitz transform of $H^0(\R^n)\equiv L^2(\R^n)$ into itself.

\medskip
\noindent
2. Let $0<s\le 1$. Then $| \cdot |$ is a continuous transform of $H^s(\R^n)$ into itself,
by general results about Nemytskii operators in Sobolev/Besov spaces, see
 {\cite[Theorem 5.5.2/3]{RS}}. 
Also it is obvious that for $u\in H^1(\R^n)$
$$
\langle -\Delta  |u|,|u|\rangle = \langle-\Delta u, u\rangle=\irn|\nabla u|^2~\!dx~,\qquad
\langle -\Delta  u^+,u^-\rangle=\irn\nabla u^+\cdot\nabla u^-~\!dx=0~\!.
$$
Here and elsewhere $u^{\pm}=\max\{\pm u,0\}=\frac{1}{2}(|u|\pm u)$, so that $u=u^+-u^-$, $|u|=u^++u^-$.
On the other hand, for $s\in(0,1)$ and $u\in H^s(\R^n)$ formula (\ref{eq:hypersing}) gives
\begin{equation}
\label{eq:identity2}
\langle \Ds  u^+,u^-\rangle= -\,C_{n,s}\iirn \frac{u^+(x)u^-(y)}{|x-y|^{n+2s}}~\!dxdy.
\end{equation}
From (\ref{eq:identity2}) we infer by the polarization identity
$$
4\langle\Ds u^+, u^-\rangle = \langle \Ds  |u|,|u|\rangle - \langle\Ds u, u\rangle
$$
that if $u$ changes sign then
\begin{equation}
\label{eq:strict0}
\langle \Ds  |u|,|u|\rangle < \langle\Ds u, u\rangle, \qquad s\in(0,1).
\end{equation}
We mention also  \cite[Theorem 6]{MNHS} for
a different proof and explanation of (\ref{eq:strict0}), that includes the case when 
$\Ds$ is replaced by the {\em Navier} (or 
{\em spectral Dirichlet}) {\em Laplacian} on
a bounded Lipschitz domain $\Omega\subset \R^n$.

\medskip
\noindent
3. 
{Let $1<s<\frac32$. The results in \cite{BM} and \cite{O} (see also 
Section 4 of the exhaustive survey \cite{RSsur}) imply
that  $| \cdot |$ is a bounded transform of $H^s(\R^n)$ into itself. That is, there
exists a constant $c(n,s)$ such that 
$$
\langle \Ds  |u|,|u|\rangle \le c(n,s) \langle\Ds u, u\rangle,\qquad u\in  H^s(\R^n)~\!.
$$
In particular, $|\cdot|$ is continuous at $0\in H^s(\R^n)$. 

It is easy to show that the assumption $s<\frac32$ can not be improved, see
 Example \ref{Ex:ex} below and \cite[Proposition p. 357]{BM}, where a more general setting
involving Besov spaces $B^{s,q}_p(\R^n)$, $s\ge 1+\frac1p$, is considered.} 

\bigskip

At our knowledge, the continuity of $|\cdot|: H^s(\R^n)\to H^s(\R^n)$, $s\in(1,\frac32)$, 
is an open problem. We can only point out the next simple result.

\begin{Proposition}
\label{P:continuity}
Let $0<\tau<s<\frac32$. Then $|\cdot|:H^s(\R^n)\to H^\tau(\R^n)$ is continuous.
\end{Proposition}

\proof
Recall that $H^s(\R^n)\hookrightarrow H^\tau(\R^n)$ for $0<\tau<s$. Actually, the Cauchy-Bunyakovsky-Schwarz inequality 
readily gives the  well known interpolation inequality
$$
\langle\Dtau v,v\rangle=\irn|\xi|^{2\tau}|{\mathcal F}[v]|^2 d\xi
\le \Big(\langle\Ds v,v\rangle\Big)^{\!\!\frac{\tau}{s}}\Big(\irn|v|^2~\!dx\Big)^{\!\!\frac{s-\tau}{s}}\!\!\!\!,\quad v\in H^s(\R^n).
$$
Since $|\cdot |$ is continuous $L^2(\R^n)\to L^2(\R^n)$
and bounded $H^s(\R^n)\to H^s(\R^n)$, the statement follows immediately.
\QED

Now we formulate our main result. 
It provides the complete proof of \cite[Theorem 1]{N} for $s$ 
below the threshold $\frac32$ and gives a positive answer to a question raised in 
\cite[Remark 4.2]{AJS} by Nicola Abatangelo, Sven Jahros and Albero Salda\~na.

\begin{Theorem}
\label{T:main}
Let $s\in(1,\frac32)$ and $u \in H^s(\R^n)$. Then formula (\ref{eq:identity2}) holds.
In particular, if $u$ changes sign then
$$
\langle \Ds  |u|,|u|\rangle > \langle\Ds u, u\rangle~\!.
$$
\end{Theorem}

Our proof is deeply based on the continuity result in Proposition \ref{P:continuity}.
The knowledge of continuity of $|\cdot|: H^s(\R^n)\to H^s(\R^n)$ could considerably simplify it.

We denote by $c$ any positive constant whose value is not important for our purposeses. Its value may
change line to line. The dependance of $c$ on certain parameters is shown in parentheses.

\section{Preliminary results and proof of Theorem \ref{T:main}}

We begin with a simple but crucial identity that has been independently pointed out
in \cite[Lemma 1]{N} and \cite[Lemma 3.11]{AJS} (without exact value of the constant).  
Notice that it holds for general fractional orders $s>0$. 

\begin{Theorem}
\label{P:Step1} Let $s>0$, $s\notin\mathbb N$. Assume that $v,w\in H^s(\R^n)$ have compact and disjoint supports. Then
\begin{equation}
\label{eq:general}
\langle\left(-\Delta_{\R^n}\!\right)^{\!s}\! {v},{w}\rangle= -\,C_{n,s}~\! 
\iirn\frac{v(x)w(y)}{|x-y|^{n+2s}}~\!dxdy.
\end{equation}
\end{Theorem}

\proof
Let $\rho_h$ be a sequence of mollifiers, and put ${w}_h:={w}*\rho_h$. Formula (\ref{eq:hypersing}) gives
\begin{multline*}
\langle\left(-\Delta_{\R^n}\!\right)^{\!s}\! {v},{w_h}\rangle=
\langle\left(-\Delta_{\R^n}\!\right)^{s-\lfloor s\rfloor}\! {v},{\left(-\Delta\right)^{\!\lfloor s\rfloor}\!w_h}\rangle\\
=\frac {C_{n,s-\lfloor s\rfloor}}2 \iirn\frac{\big(v(x)-v(y)\big)
\big(\left(-\Delta\right)^{\!\lfloor s\rfloor}\!w_h(x)-\left(-\Delta\right)^{\!\lfloor s\rfloor}\!w_h(y)\big)}
{|x-y|^{n+2(s-\lfloor s\rfloor)}}~\!dxdy.
\end{multline*}
Since for large $h$ the supports of $v$ and $w_h$ are separated, we have
$$
\langle\left(-\Delta_{\R^n}\!\right)^{\!s}\! {v},{w_h}\rangle=
-\,C_{n,s-\lfloor s\rfloor} \iirn \frac{v(x)\,\left(-\Delta\right)^{\!\lfloor s\rfloor}\!w_h(y)}
{|x-y|^{n+2(s-\lfloor s\rfloor)}}~\!dydx.
$$
Here we can integrate by parts. Using (\ref{eq:Cns}) one computes for $a>0$
$$
\Delta\, \frac {C_{n,a}}{|x-y|^{n+2a}}=\frac {C_{n,a}(n+2a)(2a+2)}{|x-y|^{n+2a+2}}=-\,\frac {C_{n,a+1}}{|x-y|^{n+2(a+1)}}
$$ 
and obtains (\ref{eq:general}) with $w_h$ instead of $w$.

Since the supports of ${v}$ and ${w}$ are separated, it is easy to pass to the limit as $h\to\infty$ and to conclude
the proof.
\QED

\begin{Remark}
\label{R:identity} Motivated by (\ref{eq:general}) and (\ref{eq:sign}),  A.I. Nazarov conjectured in \cite{N} that 
$$
\begin{array}{lll}
\langle \left(-\Delta_{\R^n}\!\right)^{\!s}\! |u|,|u|\rangle - \langle\left(-\Delta_{\R^n}\!\right)^{\!s}\! u, u\rangle<0 
& \text{ if ~~$\lfloor s\rfloor$}&\text{ is~ even;}
\\
\langle \left(-\Delta_{\R^n}\!\right)^{\!s}\! |u|,|u|\rangle - \langle \left(-\Delta_{\R^n}\!\right)^{\!s}\! u, u\rangle>0 
& \text{ if ~~$\lfloor s\rfloor$}&\text{ is~ odd}
\end{array}
$$
for any not integer exponent $s>0$ and for any changing sign function $u\in H^s(\R^n)$ such that $u^\pm \in H^s(\R^n)$.
\end{Remark}

\begin{Lemma}
\label{L:Step2}
Let $s\in(1,\frac 32)$ and $\eps>0$. If a function $u\in H^s(\R^n)$ has compact support then $(u-\eps)^+\in H^s(\R^n)$,
and 
$$
\langle \Ds (u-\eps)^+,(u-\eps)^+\rangle \le c(n,s) \langle\Ds u, u\rangle+c(n,s,{\rm supp}(u))\eps^2~\!.
$$
\end{Lemma}

\proof
Take a nonnegative function $\eta\in {\cal C}^\infty_0(\R^n)$ such that 
$\eta\equiv 1$ on ${\rm supp}(u)$. Clearly $u-\eps\eta\in H^s(\R^n)$. Hence, by Item 3 in the Introduction
we  have that 
$(u-\eps\eta)^+=(u-\eps)^+\in H^s(\R^n)$ and
\begin{eqnarray*}
\langle \Ds (u-\eps)^+,(u-\eps)^+\rangle &\le& c(n,s)~\!\langle \Ds (u-\eps\eta),u-\eps\eta\rangle\\
&\le& c(n,s)~\!\big(\langle \Ds u,u\rangle+\eps^2\langle \Ds\eta,\eta\rangle\big).
\end{eqnarray*}
The proof is complete.
\QED

In order to simplify notation, for $u:\R^n\to \R$ and $s>0$ we put 
$$\Phi^s_u(x,y)=\frac{u^+(x)u^-(y)}{|x-y|^{n+2s}}~\!.$$

\begin{Lemma}
\label{L:Step2bis} 
Let $s\in(1,\frac 32)$ and $u\in H^s(\R^n)\cap {\cal C}^0_0(\R^n)$. Then (\ref{eq:identity2}) holds, and in
particular
$\displaystyle{\Phi^s_u\in L^1(\R^n\times\R^n)}$.
\end{Lemma}

\proof
Thanks to Lemma \ref{L:Step2} we have that 
$(u^--\eps)^+ \in H^s(\R^n)\cap {\cal C}^0_0(\R^n)$ for any $\eps>0$. Next, the supports of the functions $u^+$ and $(u^--\eps)^+$ are compact and disjoint.
Thus we can apply Theorem \ref{P:Step1} to get
\begin{equation}
\label{eq:identity-eps}
\langle\Ds u^+,(u^--\eps)^+\rangle= -\,C_{n,s}\iirn \frac{u^+(x)(u(y)^--\eps)^+}{|x-y|^{n+2s}}~\!dxdy.
\end{equation}
Take a decreasing sequence $\eps \searrow 0$. 
From Lemma \ref{L:Step2}
we infer that $(u^--\eps)^+\to u^-$ weakly in $H^s(\R^n)$, as $(u^--\eps)^+\to u^-$ in $L^2(\R^n)$. Hence 
the duality product in (\ref{eq:identity-eps}) converges to the the duality product in (\ref{eq:identity2}).
Next, the integrand in the right-hand side of (\ref{eq:identity-eps}) increases to $\Phi^s_u$ a.e. on $\R^n\times\R^n$. 
By the monotone convergence theorem we get the convergence of the integrals, and the conclusion follows immediately.
\QED

\begin{Lemma}
\label{L:Step3}
Let $s\in(1,\frac 32)$ and $u\in H^s(\R^n)$. Then $\displaystyle{\Phi^s_u\in L^1(\R^n\times\R^n)}$.
\end{Lemma}

\proof
Take a sequence of functions $u_h\in {\cal C}^\infty_0(\R^n)$ such that $u_h\to u$
in $H^s(\R^n)$ and almost everywhere. Since  $\displaystyle{\Phi^s_{u_h}\to \Phi^s_u}$ a.e. on $\R^n\times\R^n$, 
Fatou's Lemma, Lemma \ref{L:Step2bis} for $u_h$ and the boundeness  of $v\mapsto v^\pm$ in $H^s(\R^n)$ give
\begin{eqnarray*}
\iirn \Phi^s_u(x,y)~\!dxdy&\le&
\liminf_{h\to \infty} \iirn \Phi^s_{u_h}(x,y)~\!dxdy=
c(n,s)\liminf_{h\to \infty}\langle\Ds u_h^+,u_h^-\rangle\\
&\le&
c(n,s)\lim_{h\to\infty} \langle\Ds u_h,u_h\rangle= c(n,s)\langle\Ds u,u\rangle,
\end{eqnarray*}
that concludes the proof.
\QED

\bigskip

\noindent
{\bf  Proof of Theorem \ref{T:main}.}
Take a sequence $u_h\in {\cal C}^\infty_0(\R^n)$ such that
$u_h\to u$ in $H^s(\R^n)$ and almost everywhere. Consider the nonnegative functions
$$
v_h:=u_h^+\wedge u^+=u^+-(u^+-u^+_h)^+~,\quad 
w_h:=u_h^-\wedge u^-=u^--(u^--u^-_h)^+.
$$
Then $v_h, w_h\in H^s(\R^n)$. Next, 
take any exponent $\tau\in(1,s)$. By Proposition \ref{P:continuity}
we have that $u^\pm-u^\pm_h\to 0$ in $H^\tau(\R^n)$; hence
$(u^\pm-u^\pm_h)^+\to 0$ in $H^\tau(\R^n)$ by Item 3 in the Introduction. Thus,
\begin{equation}
\label{eq:convergence}
v_h\to u^+~,\quad w_h\to u^-\quad \text{in $H^\tau(\R^n)$ and almost everywhere, \ as \ $h\to \infty$.}
\end{equation}

Now we take a small $\eps>0$. Recall that 
 $(v_h-\eps)^+ \in H^\tau(\R^n)$ by Lemma \ref{L:Step2}. Moreover, 
 from $0\le v_h \le u_h^+$, $0\le w_h\le u_h^-$ it follows that 
 $$
 \text{supp}((v_h-\eps)^+)\subseteq\{u_h\ge \eps\};\qquad
\text{supp}(w_h)\subseteq \text{supp}(u^-_h).
$$
In particular, 
the functions
 $(v_h-\eps)^+, w_h$ have compact and disjoint supports. Thus we can apply Theorem \ref{P:Step1} to infer
$$
\langle\Dtau (v_h-\eps)^+,{w_h}\rangle= -\,C_{n,\tau}~\! 
\iirn\frac{(v_h(x)-\eps)^+w_h(y)}{|x-y|^{n+2\tau}}~\!dxdy.
$$
We first take the limit  as $\eps\searrow 0$. The  argument in the proof of
Lemma \ref{L:Step2bis} gives
\begin{equation}
\label{eq:identity_h}
\langle\Dtau v_h,{w_h}\rangle= -\,C_{n,\tau}~\! 
\iirn\frac{v_h(x)w_h(y)}{|x-y|^{n+2\tau}}~\!dxdy.
\end{equation}
Next we push $h\to \infty$. By (\ref{eq:convergence}) we get
$$
\lim_{h\to\infty}\langle\Dtau v_h,{w_h}\rangle=\langle\Dtau u^+,u^-\rangle~\!.
$$
Further, since the integrand in the right-hand side of (\ref{eq:identity_h}) does not exceed 
$\Phi^\tau_u(x,y)$, Lemma \ref{L:Step3}, (\ref{eq:convergence}) and Lebesgue's theorem give
$$
\lim_{h\to\infty}\,\iirn\frac{v_h(x)w_h(y)}{|x-y|^{n+2\tau}}~\!dxdy=
\iirn\Phi^\tau_u(x,y)~\!dxdy~\!.
$$
Thus, we proved (\ref{eq:identity2}) with $s$ replaced by $\tau$. 
It remains to pass to the limit as $\tau\nearrow s$.
By Lebesgue's theorem, we have
\begin{eqnarray*}
\lim_{\tau\nearrow s}\,\langle \Dtau u^+,u^-\rangle &=& \lim_{\tau\nearrow s}\,\irn|\xi|^{2\tau}{\mathcal F}[u^+]~\!\overline{{\mathcal F}[u^-]}~\!d\xi \\
&=& \irn|\xi|^{2s}{\mathcal F}[u^+]~\!\overline{{\mathcal F}[u^-]}~\!d\xi=\langle \Ds u^+,u^-\rangle.
\end{eqnarray*}
Now we fix $\tau_0\in(1,s)$ and notice that $0\le \Phi^\tau_u \le \max\{\Phi^{\tau_0}_u,\Phi^s_u\}$ for any $\tau\in(\tau_0,s)$.
Therefore, Lemma \ref{L:Step3} and Lebesgue's theorem give
$$
\lim_{\tau\nearrow s}\,\iirn\Phi^\tau_u(x,y)~\!dxdy=\iirn\Phi^s_u(x,y)~\!dxdy~\!.
$$
The proof of (\ref{eq:identity2}) 
is complete. The last statement follows immediately from (\ref{eq:identity2}), polarization identity and (\ref{eq:sign}).
\QED

\begin{Example}
\label{Ex:ex}{\rm 
It is easy to construct a function $u\in {\cal C}^\infty_0(\R^n)$ such that $u^+\in H^s(\R^n)$ if and only if 
$s<\frac32$.

Take $\f\in {\cal C}^\infty_0(\R)$ satisfying 
$\f(0)=0, \f'(0)> 0$ and $x\f(x)\ge 0$ on $\R$. By direct computation one checks that $\f^+=\chi_{(0,\infty)}\f\in H^s(\R)$
if and only if $s<\frac32$. If $n=1$ we are done. If $n\ge 2$ we take
$u(x_1,x_2,\dots, x_n)=\f(x_1)\f(x_2)\dots\f(x_n)$.
}\end{Example}

\medskip
\noindent
{\bf Acknowledgements.} The first author wishes to thank Universit\'e Libre de Bruxelles for the  hospitality in
February 2016. She is grateful to Denis Bonheure, Nicola Abatangelo, Sven Jahros and Albero Salda\~na for valuable discussion
on this subject.

\end{document}